\documentclass[10pt,a4paper]{amsart}
\usepackage{amssymb,amsmath}
\newtheorem{thm}{Theorem}[section]

\newtheorem{cor}[thm]{Corollary}
\newtheorem{lem}[thm]{Lemma}

\newtheorem{rek}[thm]{Remark}

\newcommand\be{\begin{equation}}
\newcommand\ee{\end{equation}}
\newcommand\bea{\begin{eqnarray}}
\newcommand\eea{\end{eqnarray}}
\newcommand\bi{\begin{itemize}}
\newcommand\ei{\end{itemize}}
\newcommand\ben{\begin{enumerate}}
\newcommand\een{\end{enumerate}}
\newcommand\bc{\begin{center}}
\newcommand\ec{\end{center}}
\newcommand\ba{\begin{array}}
\newcommand\ea{\end{array}}

\newcommand{\R}{\ensuremath{\mathbb{R}}}

\newcommand{\foh}{\frac{1}{2}}  

\newcommand{\gd}{\delta}
\newcommand{\gD}{\Delta}
\newcommand{\gO}{\Omega}
\newcommand{\dgO}{\partial\Omega}
\newcommand{\gG}{\Gamma}

\begin{document}

\begin{center}
{\large\textbf{A UNIQUENESS THEOREM FOR THERMOACOUSTIC TOMOGRAPHY IN THE CASE OF LIMITED BOUNDARY DATA}\\
\bigskip
\medskip
Dustin \textsc{Steinhauer}\\
\bigskip
University of California, Los Angeles}
\end{center}
\bigskip
\bigskip

\textbf{Abstract.} We prove a uniqueness theorem for compactly supported initial data for the variable
speed wave equation arising in models of thermoacoustic tomography, given measurements on a part of the boundary.
The proof is based on domain of dependence arguments and D. Tataru's unique continuation theorem.

\section{Introduction}

\textit{Thermoacoustic tomography} (TAT) is a recently developed technique in medical
diagnostics.  The idea is to combine the high contrast of electromagnetic radiation with
the high resolution of ultrasound.  In utilizing this combination one hopes to overcome
the weaknesses of both imaging methods.

In TAT a brief duration microwave pulse irradiates the object.  (The practice of
\textit{photoacoustic tomography} uses a laser instead of a MW pulse, and then proceeds
similarly.)  Although the irradiation is assumed to be uniform, more energy may be
absorbed in some locations than in others.  For example, cancerous cells absorb several
times more energy than healthy cells.  For this reason, knowing the energy absorption
$f(x)$ at each point is of importance as a diagnostic tool.  The goal of TAT is to
determine $f(x)$ by using transducers on the surface of the object to measure the
pressure variation $u(x,t)$ generated by the absorption of energy.  Here $u(x,t)$ solves
the wave equation

\bea Pu = u_{tt}-c^2(x)\gD u &=& 0 \nonumber\\
u(x,0) &=&f(x)\nonumber\\
u_t(x,0) &=&0\;.\label{main}\eea
Mathematically, the problem is to recover $f(x)$ from knowledge of $u(x,t)$ on the
boundary (or part of the boundary) of the region inside the transducers.

The acoustic speed $c(x)$ is assumed to be known and smooth, having been predetermined
from ultrasound experiments.  In the case of constant $c$ there are several methods known
for reconstructing $f$, at least in the case of complete boundary data (see \cite{KuKy}
for a survey of results).  When $c$ is variable, less is known.  A reconstruction formula
was obtained in \cite{AgKu} using eigenfunctions of the Laplacian in the region $\gO$
inside the transducers.

We are interested in the case of variable $c$ and limited (partial) boundary data.  We
think it is unlikely that an exact reconstruction formula for $f$ can be obtained for general
measurement surfaces; however with knowledge of $u(x,t)$ on a sufficiently large subset
of the boundary one can conclude the uniqueness of $f$.
%

\section{Domains of Dependence}

Let $\gO$ be a bounded region in $\R^n$ with smooth boundary.  Let $c(x)$ be a positive
smooth function on $\R^n$ such that $\frac{1}{M}<c(x)<M$ for some number $M>1$.  Let $f\in
C^\infty(\R^n)$ have support inside $\gO$.

The function $c(x)$ induces a Riemannian metric on $\R^n$ by setting
$g_{ij}(x)=c^{-2}(x)\delta_{ij}$.  This metric differs from the Euclidean metric on regions where $c(x)$ is not identically equal to 1.  (In many models of TAT it is assumed that $c(x)=1$ outside a compact region; we did not find this assumption to be necessary for our purposes, so we have omitted it.)  The length of a curve $r:[a,b]\rightarrow\R^n$ in this metric is

\be L(r)=\int_a^b \frac{|r'(t)|}{c(r(t))}\;dt\;. \ee
This length depends only on the image of $r$ in $\R^n$ (i.e. is independent of
parametrization).  The \textit{geodesic distance} from $x$ to $y$ is defined to be

\be\label{distance} d(x,y) = \inf_r L(r)\;, \ee
the infimum being taken over all curves from $x$ to $y$.

If $S$ is a submanifold of $\R^n$ define the induced geodesic distance $d_S(x,y)$ for
$x,y\in S$ by taking the infimum only over curves that lie in $S$.  Note that $d(x,y)\leq
d_S(x,y)$.  We then have the following uniqueness theorem:

\begin{thm}
\label{Uniqueness}

Let $f\in C^\infty(\R^n)$ have support inside $\gO$.  Let $\gG\subset\dgO$ and suppose
$\gG$ satisfies:\\

(P) For any $x\in\gO$ there exists $p\in\gG$ such that
$d(x,p)<d_{\R^n\backslash\gO}(p,\dgO\backslash \gG)$.\\

Suppose that $T$ is large enough that $T>d(x,\gG)$ for any $x\in\gO$.  Then if $u$
satisfies

\bea u_{tt} - c^2(x)\gD u &=&0 \;\;\mbox{on}\;\; \R^n\times\R^+\\
u|_{t=0}&=&f,\; u_t|_{t=0}=0 \;.\eea
If in addition

\be u|_{\gG\times[0,T]} f=0 \ee
we have $f\equiv 0$ (and hence $u\equiv 0$).\\

\end{thm}

Note that a solution of (4)-(6) extends to all of $\R^n_x\times\R_t$ by making $u$ even in $t$.  We then have $u|_{\gG\times[-T,T]}=0$.\\

The demonstration of this has two steps.  First, we will apply Green's Theorem to show
that $u$ must vanish on the domain of dependence in $[0,T]\times\R^n\backslash\gO$ determined by $\gG\times[0,T]$.  Second,
Tataru's Uniqueness Theorem will allow us to conclude that $u$ vanishes inside $\gO$
also, and hence $f\equiv 0$.\\

\proof: Consider the following exterior initial boundary value problem:

\bea
\label{eibvp}
v_{tt} - \gD v &=&0 \;\;\mbox{on}\;\; (\R^n\backslash\gO)\times[0,T]\nonumber\\
v|_{t=0}&=&v_t|_{t=0}=0\nonumber\\
v|_{[0,T]\times\dgO}&=& g \eea
where $g$ is smooth (the solution $v(x,t)$ exists for $t\leq T$ and is unique).  Suppose
in addition that $g\equiv 0$ on $[0,T]\times\gG$.  If $g$ is taken to be equal to the
restriction of $u$ solving (4)-(6) to $\partial\gO\times[0,T]$, then $v=u$ on $(\R^n\backslash\gO)\times[0,T]$.\\

Fix a point $p\in\R^n\backslash\gO$.  We need the following fact about the function $d_{\R^n\backslash\gO} (x,p)$:

\begin{lem}
\label{Lip}

$d_{R^n\backslash\gO} (x,p)$ is a Lipschitz function of $x$ on $\R^n\backslash\gO$, and $| \nabla d_{\R^n\backslash\gO} (x,p)|\leq c^{-1}(x)$ a.e.

\end{lem}

This lemma is proved in \cite{FinPaRak} in the case $c\equiv 1$, and we extend their proof.

\begin{proof}

It is immediate from the definition in (\ref{distance}) that

\be |d_{\R^n\backslash\gO} (x,p)-d_{\R^n\backslash\gO} (y,p)|\leq d_{\R^n\backslash\gO} (x,y)\;,\ee
so it suffices to show that for every $q\in\R^n\backslash\gO$ there is a small ball $B(q,\rho)$ (defined with respect to the Euclidean metric) such that $d_{\R^n\backslash\gO} (x,y)\leq C|x-y|$ for all $x,y\in B(q,\rho)$.

First, if $q\notin\partial\gO$, let $\rho$ be small enough that $B(q,\rho)$ does not intersect $\partial\gO$.  Then if $r$ is the Euclidean segment from $x$ to $y$ in $B(q,\rho)$, we have from (\ref{distance})

\be L(r)\leq|x-y|\sup_{z\in B(q,\rho)} c^{-1}(z)\;.\ee

If $q\in\partial\gO$, then for $\rho$ small enough we can assume without loss of generality there exists a smooth function $\phi=\phi(x_1,\ldots,x_{n-1})$ such that

\be \gO\cap B(q,\rho) = \{x=(x_1,\ldots,x_n)\:|\:x_n>\phi(x_1,\ldots,x_{n-1})\}\;.\ee
As before let $r$ be the Euclidean segment from $x$ to $y$.  If $r$ does not intersect $\gO$, the same analysis applies, and $d_{\R^n\backslash\gO} (x,y)\leq C|x-y|$.  If $r$ has a portion entering $\overline{\gO}$ for the first time at $a$ and leaving it for the last time at $b$, then

\be d_{\R^n\backslash\gO} (x,y)\leq C|x-a|+C|b-y|+d_{\R^n\backslash\gO} (a,b)\leq 2C|x-y|+d_{\R^n\backslash\gO}(a,b)\;.\ee
Let $\gamma$ be the projection of the segment from $a$ to $b$ onto $\partial\gO$ parametrized by $s,\;0\leq s\leq 1$.  Then

\bea \gamma(s) &=& (1-s)(a_1,\ldots,a_{n-1},0)+s(b_1,\ldots,b_{n-1},0)\nonumber\\
&+& (0,\ldots,0,\phi((1-s)a_1+sb_1,\ldots,(1-s)a_{n-1}+sb_{n-1})\;.\eea
Hence

\be \left|\frac{d\gamma}{ds}(s)\right|\leq (n-1)|b-a|+(n-1)\sup_{x\in \partial\gO}|\nabla \phi||b-a|\leq C|b-a|\;.\ee
We conclude from this that $d_{\R^n\backslash\gO} (x,p)$ is a Lipschitz function.  By Rademacher's Theorem (see \cite{EvGa}), $d_{\R^n\backslash\gO} (x,p)$ is differentiable almost everywhere.\\

To estimate the gradient almost everywhere, let $x\in\R^n\backslash\overline{\gO}$ be a point where the gradient exists.  Fix $\gd>0$ and let $\epsilon$ be so small that $B(x,\epsilon)\subset\R^n\backslash\overline{\gO}$, and also $\sup_{|z-x|<\epsilon}c^{-1}(z)\leq c^{-1}(x)+\gd$.  Suppose $|x-y|<\epsilon$ and let $r$ parametrize the (Euclidean) segment from $x$ to $y$ by the interval $[a,b]$.  Then

\bea |d_{\R^n\backslash\gO} (x,p)-d_{\R^n\backslash\gO} (y,p)|&\leq& |d_{\R^n\backslash\gO} (x,y)|\nonumber\\
&\leq&\int_a^b\frac{|r'(t)|}{c(r(t))}\:dt\nonumber\\
&\leq&|x-y|\sup_{|z-x|<\epsilon}c^{-1}(z)\nonumber\\
&\leq&|x-y|\left(c^{-1}(x)+\gd\right)\;.\eea
Since this holds for every $\gd>0$, we must have $\nabla d_{\R^n\backslash\gO} (x,p)\leq c^{-1}(x)$ a.e.

\end{proof}

Continuing with the proof of Theorem \ref{Uniqueness}, let $x_0\in\gO$. By Property (P) there exists $y_0\in\partial\gO$ such that $d(x_0,y_0)<d_{\R^n\backslash\partial\gO} (y_0,\gG)<T$.  Let $d_{\R^n\backslash\partial\gO} (y_0,\gG)=H$ and select a point $p\in\R^n\backslash\overline{\gO}$ close enough to $y_0$ such that the intersection of the set

\be U=\{(x,t)\in\R^n\backslash\gO\times[0,T]\:|\:d_{\R^n\backslash\gO} (x,p)+t<H\}\ee
with $\partial\gO\times[0,T]$ is contained inside $\gG\times[0,T]$.  Such a $p$ can be chosen precisely because of Property (P).  We wish to show that $v$ solving (\ref{eibvp}) vanishes in $U$.  Such a set $U$ is called a \textit{domain of dependence} (see \cite{KaKuLa}).  A slight modification of the usual Green's Theorem argument will allow us to make this conclusion.  The argument must be modified because, as the authors point out in \cite{FinPaRak}, $\partial U$ is only a Lipschitz surface.  In order to overcome this difficulty, we will approximate $U$ from the inside by regions with piecewise smooth boundary.

$d_{\R^n\backslash\gO} (x,p)$ is a Lipschitz function, so it can be approximated by a smooth function in the sense that, for any $\epsilon>0$ there exists a smooth function $f_\epsilon$ on $\R^n\backslash\gO\cap\{d_{\R^n\backslash\gO} (x,p)\leq T\}$ such that $|d_{\R^n\backslash\gO} (x,p)-f_\epsilon(x)|<\epsilon$ and $|\nabla f_\epsilon(x)|\leq c^{-1}(x)+\epsilon$ (see \cite{Az}).  We can also approximate $(1-\gd)d_{\R^n\backslash\gO} (x,p)$ by a smooth function $f_{\epsilon,\gd}$ so that $|(1-\gd)d_{\R^n\backslash\gO} (x,p)-f_{\epsilon,\gd}(x)|<\epsilon$ and $|\nabla f_{\epsilon,\gd}(x)|\leq(1-\gd) c^{-1}(x)+\epsilon$ for any $\gd>0$.  For $\epsilon$ sufficiently small relative to $\gd$, we have the important estimate

\be |\nabla f_{\epsilon,\gd}(x)|\leq c^{-1}(x)\;.\ee

We next define the domains of dependence on which we will apply Green's Theorem. For $0<h\leq H$ let $R=R_{\epsilon,\gd,h}=\{(x,t)\in\R^n\backslash\gO\times[0,h]\:|\:f_{\epsilon,\gd}(x)+t<h\}$.  Then $R$ is a domain with piecewise smooth boundary, and for $\epsilon$ sufficiently small, we still have $\partial R\cap\partial\gO\subset\gG$.  In fact, $\partial R$ consists of three pieces (Figure 1):

\begin{itemize}

\item $\Sigma_1$: A slice of $\{t=0\}$\\

\item $\Sigma_2$: A spacelike surface which is part of $\{(x,t)\:|\: f_{\epsilon,\gd}(x)+t=h\}$\\

\item $\Sigma_3$: A slice of $\gG\times[0,h]$\\

\end{itemize}
Here, \textit{spacelike} means $c(x)|N_x|\leq|N_t|$ for all normal vectors $(N_x,N_t)$.  The proper terminology should be \textit{spacelike or null}, but we shall simply write spacelike for brevity.

Using Green's Theorem (see for example \cite{Tay}) we have

\bea 0&=&\int_R v_t(c^{-2}(x)v_{tt}-\Delta v)\:dV\:dt \nonumber\\
&=&\int_R c^{-2}(x)\frac{\partial}{\partial t}\left(\foh v_t^2\right)\:dV\:dt +\int_R \langle \nabla_xv_t,\nabla_xv\rangle\:dV\:dt\nonumber\\
&-&\int_R \mbox{div}_x(v_t\: \nabla_xv)\:dV\:dt\;\nonumber\\
&=&\foh\int_{\Sigma_1\cup\Sigma_2}(c^{-2}(x)v_t^2 + \langle \nabla_xv,\nabla_xv\rangle)\:\omega - \int_{\Sigma_2\cup\Sigma_3} v_t\frac{\partial v}{\partial \nu_x}\:dS_t\:dt\;.\eea
If $dS$ is the measure on $\partial R$ induced by $dV\:dt$, then $\omega = N_t\:dS$ and $dS_t\:dt=|N_x|\:dS$ where $(N_x,N_t)$ is the outward unit normal to $\partial R$.  As a result,

\be \label{yourmom} \int_{\Sigma_2}\left((c^{-2}(x)v_t^2+|\nabla_xv|^2)|N_t|-2v_t\frac{\partial v}{\partial \nu_x}|N_x|\right)\:dS=0\;.\ee
$v$, and hence $v_t$, vanishes on $\Sigma_3$, so

\be \int_{\Sigma_3}v_t\frac{\partial v}{\partial \nu_x}|N_x|\:dS=0\;.\ee
On the other hand, $\Sigma_2$ is spacelike, so since we have

\be 2\left|c^{-1}(x)v_t\frac{\partial v}{\partial \nu_x}\right|\leq c^{-2}(x)v_t^2+|\nabla_xv|^2\;,\ee
we obtain a positive-definite form in $(v_t,\nabla_xv)$:

\be (c^{-2}(x)v_t^2+|\nabla_xv|^2)|N_t|-2v_t\frac{\partial v}{\partial \nu_x}|N_x| \geq0 \ee
with equality only when $(v_t,\nabla_xv)\equiv0$.  Therefore (\ref{yourmom}) implies

\be \int_{\Sigma_2}\left((c^{-2}(x)v_t^2+|\nabla_xv|^2)|N_t|-2v_t\frac{\partial v}{\partial \nu_x}|N_x|\right)\:dS=0 \;.\ee
As a result, $(v_t,\nabla_xv)=0$ on $\Sigma_2$.  Letting $h$ vary from 0 to $H$, observing that the surfaces $\Sigma_2(h)$ sweep out $R_{\epsilon,\gd,H}$, and using that the initial data vanish, we conclude that $v\equiv0$ on $R_{\epsilon,\gd,H}$.

Now let $(x,t)\in U$.  Then $d_{\R^n\backslash\gO} (x,p)+t<H$, so for $\gd$ sufficiently small, \newline $(1-\gd)d_{\R^n\backslash\gO} (x,p)+t<H$.  As a result,

\be d_{\R^n\backslash\gO} (x,R_{\epsilon,\gd,H})<\epsilon \ee
so by the continuity of $v$, $v\equiv0$ on $U$.

Extending $v$ to be even in time, we also have $v\equiv0$ on $-R:=\{(t,x)|(-t,x)\in R$\}.  If we shrink $H$ by any sufficiently small amount, then there exists $\rho>0$ such that $v\equiv0$ on $B(p,\rho)\times[-H,H]$ while maintaining $d(x_0,p)<H$.  This is the fact we will require in the next section.

%
%
%
%
%
%
%

\section{Tataru's Theorem and Unique Continuation}

In this section we will make use of the following special case of a uniqueness theorem of D. Tataru:

\begin{thm}
\emph{[Tataru]}
\label{Tataru}
Let $u\in \mathcal{D}(\R,H^1(\R^n))$ be a solution to

\be \partial_t^2u-\sum_{i,j\leq1}^n a_{ij}(x)\partial_{x_i}\partial_{x_j}u=0\ee
where the coefficients $a_{ij}$ are smooth.  Suppose $S$ is a noncharacteristic
hypersurface containing a point $y$, and that $u$ vanishes on one side of $S$ near
$y$.  Then $u$ must vanish in a neighborhood of $y$.

\end{thm}

We will use Tataru's Theorem to prove the following fact about solutions which vanish in
a cylinder:

\begin{lem}
\label{Kach}

Suppose $u$ solves $Pu=0$ and $u$ vanishes on the cylinder
$B(z,\rho)\times[-D,D]$.  Here, $B(z,\rho)$ is defined with respect to the metric $g$, i.e. $B(z,\rho)=\{x\in\R^n\:|\: d(x,z)<\rho\}$.  Suppose $D$ is less than the injectivity
radius of the metric $g_{ij}$ at $z$.  Then $u$ vanishes on $X=\{(x,t)\:|\:
d(x,z)+|t|<D\}$ (see Figure 2).

\end{lem}

\begin{proof}

Let $\chi=\chi_\epsilon\in C^\infty_0(\R)$ be supported in
$[-\epsilon,\epsilon]$, satisfy $\chi(t)=1$ for $|t|\leq \frac{\epsilon}{2}$, have
$\chi(t)\leq\foh$ for $|r|\geq\frac{3}{4}\epsilon$ and have $|\chi'(t)|\leq
\frac{3}{\epsilon}$ for all $t$.  Define for $0<r<D$ and $\epsilon$ small:

\be X_{r,\epsilon}=\{(x,t)\:|\: d(x,z)\frac{D}{r} +
|t|(1-\chi(t))+\frac{\epsilon}{2}\chi(t)<D\} \;.\ee
Note that for $r<\rho$, $X_{r_,\epsilon}\subset B(z,\rho)\times[-D,D]$, so $u$ vanishes
on $X_{r,\epsilon}$ for $r<\rho$.

$\partial X_{r,\epsilon}$ is a smooth hypersurface except at the points $(z,\pm D)$.  Let
$(x,t)\in\partial X_{r,\epsilon},\; |t|<D$.  We can compute the normal to $\partial
X_{r,\epsilon}$ at $(x,t)$ with respect to the metric $g+dt^2$ by considering the
gradient of the function

\be F(x,t) = d(x,z)\frac{D}{r} + |t|(1-\chi(t))+\frac{\epsilon}{2}\chi(t)\;.\ee
Since $x$ lies inside a system of normal coordinates about $z$, we have
$|\nabla_xF(x,t)|=\nabla_xd(x,z)\frac{D}{r}=c^{-1}(x)\frac{D}{r}$.  We can also estimate
$d_tF(x,t)$ by computing

\be d_tF(x,t)=(1-\chi(t))\frac{|t|}{t}-\chi'(t)\left(|t|-\frac{\epsilon}{2}\right)\;.\ee
If $|t|\geq\epsilon$, $|d_tF|=1$.  If $\frac{3}{4}\epsilon<|t|<\epsilon$, then
$1-\chi(t)>\foh$, and

\be |d_tF(x,t)|\leq \left|1-\chi(t)-\frac{3}{\epsilon}\cdot\frac{\epsilon}{2}\right|\leq 1;.\ee
If $\foh\epsilon<|t|\leq\frac{3}{4}\epsilon$, then
$|t|-\foh\epsilon<\frac{1}{4}\epsilon$, and

\be |d_tF(x,t)|\leq
|\chi'(t)|\left(|t|-\foh\epsilon\right)\leq\frac{3}{\epsilon}\cdot\frac{\epsilon}{4}<1\;.\ee
Obviously if $|t|\leq \foh\epsilon$, $d_tF(x,t)=0$.

The symbol of the operator $P$ is

\be p(x,\xi,\tau)=\tau^2-c^2(x)|\xi|^2\;.\ee
At the point $(x,t)$, the normal to $\partial X_{r,\epsilon}$ is $(x,t,\nabla_xF,d_tF)$
and satisfies

\be p( x,t,\nabla_xF,d_tF) \leq 1-c^2(x)\left(c^{-2}(x)\frac{D^2}{r^2}\right)<0 \ee
so $\partial X_{r,\epsilon}$ is noncharacteristic at $(x,t)$.

Let $R = \sup_r \;\{r\:\mid\: u|_{X_{r,\epsilon}}=0\;\; \forall \epsilon>0\}$.  By
continuity $u=0$ on $X_{R,\epsilon} \;\forall \epsilon>0$.  But if $R<D$, the boundary of
$X_{R,\epsilon}$ is noncharacteristic, so $u$ must vanish in a neighborhood of
$X_{R,\epsilon}$; hence $R=D$.  The regions $X_{r,\epsilon}$ exhaust $X$, completing the
proof.

\end{proof}

A proof of this lemma can also be found in \cite{KaKuLa}.\\

We finish the proof of Theorem \ref{Uniqueness} as follows: The function $c(x)$ is
bounded away from 0, and all of its first and second derivatives are bounded functions. 
This implies $g_{ij}$ is a metric whose curvature components are bounded, and as a
result, there exists a number $\gd>0$ such that $\forall x\in\R^n$ there are geodesic
normal coordinates about $x$ covering a region $\{y\:|\:d(x,y)<\gd\}$.  We obtain the following corollary of Lemma \ref{Kach}:

\begin{cor}

Let $r>0$.  If $u=0$ on $B(p,r)\times[-H, H]$, then $u=0$ on the set

\be Y = \{(x,t)\:|\: d(x,z)<r+\operatorname{min}(H-|t|,\gd)\}\;.\ee

\end{cor}

\begin{proof}

Let $(x,t)\in Y$.  In the case when $t\in (-H,-H+\gd)$, we apply Lemma \ref{Kach} knowing $u$ vanishes on $B(p,r)\times [-H,-H+2\gd]$ to conclude $u(x,t)=0$.  Similarly, if $t\in (H-\gd,H)$ we apply the lemma knowing $u$ vanishes on $B(p,r)\times [H-2\gd,H]$.  Finally, if $-H+\gd\leq t\leq H-\gd$, we find $u(x,t)=0$ by applying the lemma, knowing $u$ vanishes on $B(p,r)\times [t-\gd, t+\gd]$.\end{proof}


Repeated applications of this corollary to the cylinders

\be B(p,\rho+n\gd)\times [-(H-n\gd),H-n\gd] \nonumber\ee
show that if $u$ vanishes on $B(p,\rho)\times[-H,H]$ for some $\rho>0$, then $u$ vanishes at $(x,t)$ if $d(x,z)+|t|\leq H$.  In particular, $f(x_0)=u(x_0,0)=0$.  Hence $f\equiv 0$, concluding the proof of Theorem \ref{Uniqueness}.$\;\;\;\;\;\;\;\;\;\;\;\;\;\;\;\;\;\;\;\;\;\;\;\;\;\;\;\;\;\;\;\;\;\;\;\;\;\;\;\;\;\;\;\;\;\;\;\;\;\;\;\;\;\;\;\;\;\;\;\;\;\;\;\;\;\;\;\;\;\;\;\;\;\;\;\;\;\;\;\;\;\;\;\;\;\;\;\;\;\;\;\;\;\;\;\;\;\;\Box$

\begin{rek}
\begin{center}
This work was done several months ago and was included my preliminary Ph.D. examination in June, 2008.  It appears my results are similar to the results contained in a paper of P. Stefanov and G. Uhlmann
(http://arxiv.org/PS\_cache/arxiv/pdf/0902/0902.1973v1.pdf)

that was just announced.\end{center}
\end{rek}

\bigskip
\bigskip

\noindent
\textsc{UCLA Department of Mathematics, Los Angeles, CA 90095-1555, USA}\\
\textit{E-mail address:} \textbf{dsteinha@math.ucla.edu}


\begin{thebibliography}{widest-label}

\bibitem[AgKu]{AgKu}Agranovsky, M. \& Kuchment, P. 2007 Uniqueness of reconstruction and
an inversion procedure for thermoacoustic and photoacoustic tomography
with variable sound speed. Inverse Problems 23, 2089-2102.

\bibitem[AgQu]{AgQu}Agranovsky, M. \& Quinto, E. T. 1996 Injectivity sets for the Radon
transform over circles and complete systems of radial functions. Journal
of Functional Analysis, 139, 383-414.

\bibitem[Az]{Az}Azagra, D., Ferrera, J., L\'opez-Mesas, F. \& Rangel, Y. (2006) Smooth approximation of Lipschitz functions on Riemannian manifolds. Journal of Mathematical Analysis and Applications, Vol. 326, No. 2, pp. 1370-1378.

\bibitem[Bey]{Bey}Beylkin, G. Imaging of discontinuities in the Inverse Scattering
Problem, G. J. Math Phys. 26(1), January 1985.

\bibitem[Du]{Du} Duistermaat, J.J. Fourier Integral Operators (lecture notes).  Courant
Institute of Mathematical Sciences, 1973.

\bibitem[EvGa]{EvGa} Evans, L.C. \& Gariepy, R.F. Measure Theory and Fine Properties of Functions, CRC Press,
Boca Raton, FL, 1992.

\bibitem[FinPaRak]{FinPaRak} Finch, D., Patch, S. \& Rakesh (2004) Determining a function from its
mean values over a family of spheres, SIAM J. Math. Anal. Vol. 35 No. 5 pp. 1213-1240.

\bibitem[Ho]{Ho} H\"ormander, L. Fourier Integral Operators I.  Acta Math 127 (1971) 79-183.

\bibitem[KaKuLa]{KaKuLa} Kachalov, A., Kurylev, Y., Lassas, M. Inverse Boundary Spectral Problems, Boca Raton: Chapman and Hall, 2001.

\bibitem[KuKy]{KuKy}Kuchment, P. \& Kuyansky, L. Mathematics of Thermoacoustic
Tomography, European J. Appl. Math., v. 19, 2008, pp. 1-34.

\bibitem[Rak]{Rak} Rakesh. (1988) A Linearised inverse problem for the wave
equation, Communications in Partial Differential Equations, 13:5, 573 - 601

\bibitem[Tar]{Tar} Tataru, D. 1995 Unique continuation for solutions to pde's;
between hörmander's theorem and holmgren' theorem, Communications in Partial
Differential Equations, 20:5, 855-884.

\bibitem[Tay]{Tay} Taylor, D. Partial Differential Equations, Vol I \& II,
Springer-Verlag, 1996.

\bibitem[Xu]{Xu} Xu, Y., Wang, L., Ambartsoumian, G. \& Kuchment, P. 2004 Reconstructions
in limited view thermoacoustic tomography. Medical Physics, 31(4), 724-733.

\end{thebibliography}
\end{document}